\documentclass[12pt]{amsart}
\usepackage{amssymb, amsmath}

\newtheorem{theorem}{Theorem}[section]
\newtheorem{proposition}[theorem]{Proposition}
\newtheorem{corollary}[theorem]{Corollary}
\newtheorem{lemma}[theorem]{Lemma}

\theoremstyle{definition}
\newtheorem{definition}[theorem]{Definition}

\theoremstyle{remark}

\begin{document}

\newcommand{\fra}[1]{{\mathfrak{#1}}}
\newcommand{\Gr}{\text{Gr}}
\newcommand{\grh}{{\text{gr}}_{F}}
\newcommand{\ra}{\rightarrow}
\newcommand{\Ra}{\Rightarrow}
\newcommand{\surj}{\twoheadrightarrow}
\newcommand{\lra}{\longrightarrow}
\newcommand{\pa}{\partial}
\newcommand{\noi}{\noindent}
\newcommand{\PP}{\mathbf{P}}
\newcommand{\PPP}{\PP_{\text{sub}}}
\newcommand{\RR}{\mathbf{R}}
\newcommand{\NN}{\mathbf{N}}
\newcommand{\lef}{\mathbf{L}}
\newcommand{\hyp}{\mathbf{H}}
\newcommand{\ZZ}{\mathbf{Z}}
\newcommand{\CC}{\mathbf{C}}
\newcommand{\QQ}{\mathbf{Q}}
\newcommand{\bin}[2]{ {{#1} \choose {#2}} }
\newcommand{\OO}{\mathcal{O}}

\newcommand{\cH}{\mathcal{H}}
\newcommand{\MM}{\mathcal{M}}
\newcommand{\KK}{\mathcal{K}}
\newcommand{\II}{\mathcal{I}}
\newcommand{\LL}{\mathcal{L}}
\newcommand{\FF}{\mathcal{F}}
\newcommand{\GG}{\mathcal{G}}
\newcommand{\EE}{\mathcal{E}}
\newcommand{\cM}{\mathcal{M}}
\newcommand{\cN}{\mathcal{N}}
\newcommand{\DD}{\mathcal{D}}
\newcommand{\cA}{\mathbf{A}}
\newcommand{\ti}[1]{\tilde{#1}}
\newcommand{\ef}{\rm{\ if\  }}
\newcommand{\eps}{\epsilon}
\newcommand{\sbl}{\vskip 3pt}
\newcommand{\lbl}{\vskip 6pt}
\newcommand{\rk} {\text{rank }}
\newcommand{\osect}{\mathbf{0}}
\newcommand{\HH}[3]{H^{{#1}} \big( {#2} , {#3} \big) }
\newcommand{\hh}[3]{h^{{#1}} \big( {#2} , {#3} \big) }
\newcommand{\OP}[1]{\OO_{\PP({#1})}}
\newcommand{\fall}{ \ \ \text{ for all } \ }
\newcommand{\rndup}[1]{ \ulcorner {#1} \urcorner }
\newcommand{\rndown}[1] {\llcorner {#1} \lrcorner}
\newcommand{\hgt}{\rm{height }}
\newcommand{\MI}[1]{\mathcal{J} ( {#1} ) }
\newcommand{\MIP}[1]{\mathcal{J}_+ ( {#1} ) }
\newcommand{\ZMI}[1]{\text{Zeroes}\big( \MI{ {#1} }
          \big)}
\newcommand{\BI}[1]{  \mathfrak{b} \big( {#1} \big) }
\newcommand{\ord}{\text{ord}}
\newcommand{\codim}{\text{codim}}
\newcommand{\mult}{\text{mult}}
\newcommand{\Supp}{\text{Supp}}
\newcommand{\defi}{{\text{def}}}
\newcommand{\Ex}{\text{Ex}}
\newcommand{\QT}[2]{{#1}<{#2}>}
\newcommand{\Div}{\text{Div}}
\newcommand{\num}{ \equiv_{\text{num}} }
\newcommand{\lin}{\equiv}
\newcommand{\Qnum}{\equiv_{\text{num},\QQ}}
\newcommand{\Qlin}{\equiv_{\text{lin},\QQ}}
\newcommand{\dra}{\dashrightarrow}
\newcommand{\Bl}{\text{Bl}}
\newcommand{\lc}{\text{lc}}
\newcommand{\length}{\text{length}}
\newcommand{\lam}{\lambda}
\newcommand{\al}{\alpha}
\newcommand{\gr}{\text{gr}}
\newcommand{\te}{\text}
\newcommand{\Sp}{\text{Sp}}
\newcommand{\HSp}{\text{HSp}}
\newcommand{\db}{\underline{\Omega}\,\dot{}}

\title[On the $V$-filtration of $\DD$-modules]{On the $V$-filtration of $\DD$-modules}
\author{Nero Budur}
\address{Department of Mathematics, Johns Hopkins University,
Baltimore, MD 21218-2686, USA} \email{nbudur@math.jhu.edu}

\begin{abstract}
In this mostly expository note we give a down-to-earth
introduction to the $V$-filtration of M. Kashiwara and B.
Malgrange on $\DD$-modules. We survey some applications to
generalized Bernstein-Sato polynomials, multiplier ideals, and
monodromy of vanishing cycles.
\end{abstract}

\maketitle

The $V$ filtration on $\DD$-modules was introduced by M. Kashiwara
and B. Malgrange to construct vanishing cycles in the category of
(regular holonomic) $\DD$-modules. Our aim is to give a
down-to-earth introduction to this notion and describe some
applications. The first application is to the generalized
Bernstein-Sato polynomials introduced in \cite{BMS}. Following G.
Lyubeznik, we extend a finiteness result on the set of these
polynomials. Then we describe applications to multiplier ideals
(\cite{BS}, \cite{BMS}) and to monodromy of vanishing cycles and
Hodge spectrum (\cite{B}, \cite{BS}).

We thank M. Musta\c{t}\u{a} who wrote and provided us with
preliminary notes on the $V$ filtration. We thank M. Saito for
clarifications on many issues. Most of what I learned about the
$V$-filtration is from discussions with them and with L. Ein.

\section{Basics}

In this section we introduce the filtration $V$ and prove a few
consequences assuming its existence. For a complete account on the
$V$ filtration consult \cite{Ma}, \cite{Ka}, \cite{Lau},
\cite{Sab}.

Let $X$ be a smooth complex variety. The sheaf $\DD_X$ of
algebraic differential operators on $X$ is generated locally by
multiplication by functions and by the tangent vector fields. If
$X=\cA^n$ is the affine $n$-space, then $\DD_X$ is the Weyl
algebra  $$A_n(\CC)=\CC[x_1,\ldots,x_n,\pa_1,\ldots,\pa_n],$$
where $\pa_i=\pa/\pa_{x_i}$ and $\pa_ix_j-x_j\pa_i=\delta_{i,j}$.
In these notes we will consider only left $\DD_X$-modules, the
most important for our applications being $\OO_X$. We will
frequently work locally without specifically assuming that $X$ is
affine.

The $V$-filtration of M. Kashiwara and B. Malgrange on
$\DD_X$-modules is defined with respect to some closed subvariety
$Z\subset X$.

\medskip

{\it Case $Z\subset X$ smooth.} First we consider smooth closed
subvarieties $Z\subset X$. Let $\II\subset\OO_X$ denote the ideal
of $Z$. In local coordinates, write
$$X=\{(x,t)\},\ Z=\{x\}=\{t=0\},$$ with $x=x_1,\ldots, x_n$, and
$t=t_1,\ldots, t_r$. Then $$\DD_X=\CC[x,t,\pa_x,\pa_t].$$

The $V$-filtration on $\DD_X$ is defined by
$$V^j\DD_X=\{\ P\in\DD_X\ |\ P\II^i\subset \II^{i+j}\text{ for all
}i\in\ZZ\ \},$$ with $j\in\ZZ$ and $\II^i=\OO_X$ for $i\le 0$. Locally,
$$V^j\DD_X=\sum_{|\beta|-|\gamma|\ge j}  h_{\al,\beta,\gamma}(x)\pa_x^\al
t^\beta\pa_t^\gamma.$$ Here we use vectorial indices for
monomials, and  $|\beta|=\sum_i\beta_i$. A computation with local
coordinates shows: \medskip

(i) $V^{j_1}\DD_X\cdot V^{j_2}\DD_X\subset V^{j_1+j_2}\DD_X$, with
equality if $j_1,j_2\ge 0$;

(ii) $V^j\DD_X=\II^j\cdot V^0\DD_X\cdot\DD_{X,-j}=\DD_{X,-j}\cdot
V^0\DD_X\cdot \II^j$, where $\DD_{X,j}\subset\DD_X$ are the
operators of order $\le j$, and $\II^j=\DD_{X,j}=\OO_X$ for $j\le
0$.

\begin{definition}\label{dv}
The filtration $V$ along $Z$ on a coherent left $\DD_X$-module $M$
is an exhaustive decreasing filtration of coherent
$V^0\DD_X$-submodules $V^\al:=V^\al M$, such that:

(i) $\{V^\al\}_\al$ is indexed left-continuously and discretely by
rational numbers, i.e. $V^\al=\cap_{\beta<\al}V^\beta$, every
interval contains only finitely many $\al$ with $\Gr_V^\al\ne 0$,
and these $\al$ must be rational. Here, $\Gr^\al_V =V^\al
/V^{>\al}$, where $V^{>\al}=\cup_{\beta>\al}V^\beta $.

(ii) $t_jV^\al \subset V^{\al+1}$, and $\pa_{t_j}V^\al \subset
V^{\al-1}$ for all $\al\in \QQ$, i.e. $(V^i\DD_X)(V^\al M) \subset
 V^{\al+i}M$;

(iii) $\sum_jt_jV^\al =V^{\al+1}$ for $\al\gg 0$;

(iv) the action of $\sum_j\pa_{t_j}t_j-\al$ on $\Gr^\al_V $ is
nilpotent on $X$.
\end{definition}

All conditions are independent of the choice of local coordinates.

\begin{theorem}( M. Kashiwara, B. Malgrange )\label{existence} The filtration $V$ along $Z$ exists if $M$ is
regular holonomic and quasi-unipotent.
\end{theorem}

It is beyond our scope to introduce the theory of holonomic
systems of differential operators with regular singularities (see
\cite{Bo}, \cite{Ka2}). It suffices to say that all the
$\DD$-modules considered in the applications are regular holonomic
and quasi-unipotent.


\begin{proposition} The $V$-filtration along $Z$  is unique.
\end{proposition}
\begin{proof} Let $\widetilde{V}$ be another filtration on $M$ satisfying
Definition \ref{dv}. By symmetry it suffices to show that
$V^{\alpha}\subset \widetilde{V}^{\alpha}$ for every $\alpha$.
Suppose that $\alpha\neq\beta$ and consider
$$V^{\alpha}\cap\widetilde{V}^{\beta}/(V^{>\alpha}\cap
\widetilde{V}^{\beta})+(V^{\alpha}\cap\widetilde{V}^{>\beta}).$$
Since both filtrations satisfy \ref{dv}-(iv), it follows that both
$(\sum_j\partial_{t_j}t_j-\alpha)$ and
$(\sum_j\partial_{t_j}t_j-\beta)$ are nilpotent on this module.
Hence the module is zero.

We show now that for every $\alpha$ we have
\begin{equation}\label{formula11}
V^{\alpha}\subset V^{>\alpha}+\widetilde{V}^{\alpha}.
\end{equation}
Fix $w\in V^{\alpha}$. By exhaustion, there is $\beta\ll 0$ (in
particular $\beta<\alpha$) such that $w\in \widetilde{V}^{\beta}$.
By what we have already proved, we may write $w=w_1+w_2$, with
$w_1\in V^{>\alpha}$ and $w_2\in V^{\alpha}\cap
\widetilde{V}^{>\beta}$. If we replace $w$ by $w_2$, then the
class in $V^{\alpha}/V^{>\alpha}$ remains unchanged, but we may
choose a larger $\beta$. We can repeat the process as long as
$\beta<\alpha$. Since the $\widetilde{V}$-filtration is discrete,
we can repeat the process until we have $\beta\geq\alpha$. Hence
the class of $w$ in $V^{\alpha}/V^{>\alpha}$ can be represented by
an element in $\widetilde{V}^{\al}$, and we get
$(\ref{formula11})$.

Since the $V$-filtration is discrete, a repeated application of
$(\ref{formula11})$ shows that for every $\beta\geq\alpha$ we have
$V^{\alpha}\subset V^{\beta}+\widetilde{V}^{\alpha}$. We deduce from \ref{dv}-(iii)
that if we fix $\beta\gg 0$, then
\begin{equation}\label{fl2}
V^{\alpha}\subset \II^q\cdot
V^{\beta}+\widetilde{V}^{\alpha}
\end{equation}
for big enough $q$, where $\II\subset\OO_X$ is the ideal pf $Z$. By coherence, $V^\beta=\sum
V^0\DD_X\cdot u_i$ for finitely many $u_i$. By exhaustion, there exists  some
$\gamma\in\ZZ$ such that $\widetilde{V}^\gamma$ contains the $u_i$, hence also
$V^\beta$. By \ref{dv}-(ii), for $q$ with $q+\gamma\ge\al$ we have  $\II^q\widetilde{V}^\gamma\subset
\widetilde{V}^\al$. Thus $\II^q V^\beta\subset\widetilde{V}^\al$. Hence by (\ref{fl2}) we have $V^\al\subset\widetilde{V}^\al$.

\end{proof}

\medskip
{\it Case $Z\subset X$ arbitrary}. Now let $X$ be a smooth complex
variety and $Z\ne X$ a closed subscheme. Suppose $f_1,\ldots,
f_r\in\OO_X$ generate the ideal $\II\subset\OO_X$ of $Z$. Let
$i:X\ra X\times\cA^r=Y$ be the embedding $x\mapsto
(x,f_1(x),\ldots,f_r(x))$. Let $t_j:Y\ra\cA^1$ be the projection
with $t_j\circ i=f_j$.

Let $N$ be a $\DD_X$-module. Let $M=i_*N$, where $i_*$ is the
direct image for left $D$-modules. Working out the definition of
the direct image (e.g. \cite{Bo}), one gets
$M=N\otimes\CC[\pa_{t_1},\ldots,\pa_{t_r}]$ with the left
$\DD_Y$-action given as follows. Let $x_1,\ldots,x_n$ be local
coordinates on $X$. For $g\in\OO_X, m\in N$, and
$\pa_t^\nu=\pa_{t_1}^{\nu_1}\ldots\pa_{t_r}^{\nu_r}$,

\begin{align*}
&g(m\otimes \pa_t^\nu)=gm\otimes \pa_t^\nu,\ \ \pa_{x_i}(m\otimes
\pa_t^\nu)=\pa_{x_i}m\otimes \pa_t^\nu-\sum_j\frac{\pa f_j}{\pa
x_i}m\otimes\pa_{t_j}\pa_t^\nu\ ,\\
&\pa_{t_j}(m\otimes \pa_t^\nu)=m\otimes\pa_{t_j}\pa_t^\nu,\ \
t_j(m\otimes \pa_t^\nu)=f_jm\otimes \pa_t^\nu-\nu_jm\otimes
(\pa_t^\nu)_{j}\ ,
\end{align*}

\medskip
\noindent where $(\pa_t^\nu)_j$ is obtained from $\pa_t^\nu$ by
replacing $\nu_j$ with $\nu_j-1$. If $N$ satisfies the
requirements of Theorem \ref{existence}, then also $M$ does.

\begin{definition}\label{graph} The $V$-filtration along $Z$ on $N(=N\otimes 1)$ is defined by $V^\al N=(N\otimes 1)\cap V^\al
M,$ for $\al\in \QQ$ and $V$ on $M$ taken along $X\times \{0\}$.
\end{definition}

\begin{proposition} The definition above depends on the ideal
$\II$ of $Z$ in $X$ and not on the particular generators chosen.
\end{proposition}
\begin{proof} (cf. \cite{BMS}-2.7) Suppose $g_1,\ldots, g_{r'}\in\OO_X$ also generate $\II$, with $g_j=\sum a_{ij}f_i$, $a_{ij}\in\OO_X$. Let
$i':Y=X\times\cA^r\ra Y'=X\times\cA^r\times\cA^{r'}$ be the
embedding sending $(x,t)$ to $(x,t,t')$, where $t'_j=\sum
a_{ij}(x)t_i$, $j=1,\ldots ,r'$. The crucial fact here is that the
image of $X\times\{0\}$ is $X\times\{0\}\times\{0\}$. Working
locally, we can assume that $x,t,u$ is a local coordinate system
on $Y'$ such that $Y=\{u=0\}$, $X=\{t=u=0\}$. Hence $M'=i'_*M$ can
be written as $M\otimes\CC[\pa_{u_1},\ldots ,\pa_{u_{r'}}]$ with
left $\DD_{Y'}$-action as above. Note that some simplifications
occur: $\pa_{x_i}(m\otimes\pa_u^\nu)=\pa_{x_i}m\otimes\pa_u^\nu$,
$\pa_{t_i}(m\otimes\pa_u^\nu)=\pa_{t_i}m\otimes\pa_u^\nu$, and
$u_j(m\otimes\pa_u^\nu)=-\nu_jm\otimes(\pa_u^\nu)_j$, where $m\in
M$, $\pa_u^\nu=\pa_{u_1}^{\nu_1}\ldots\pa_{u_{r'}}^{\nu_{r'}}$,
and $(\pa_u^\nu)_j$ is obtained from $\pa_u^\nu$ by replacing
$\nu_j$ with $\nu_j-1$.

The claim follows if we show that
$$V^\al M'=\sum_{\nu\,\in\,\NN^{r'}} V^{\al+|\nu|}M\otimes\CC[\pa_u^\nu]$$
is the $V$-filtration on $M'$ along $X\times \{0\}\times\{0\}$.

Let us check the axioms for the $V$-filtration. In local
coordinates, $V^0\DD_{Y'}$ is generated over $\OO_{Y'}$ by the
$\pa_{x_i}$, and the $v\pa_{w}$ with
$v,w\in\{t_1,\ldots,t_r,u_1,\ldots, u_{r'}\}$. From definition,
these actions are well-defined on $V^\al M'$. To show that $V^\al
M'$ is coherent over $V^0\DD_{Y'}$, it is enough to show that
$V^\al M'$ is locally finitely generated since $V^0\DD_{Y'}$ is
coherent. Since $V^\al M$ is locally finitely generated over
$V^0\DD_Y$,  we have that for $c\gg 0$, $\sum_{|\nu|\le
c}V^{\al+|\nu|}M\otimes\CC[\pa_u^\nu]$ is locally finitely
generated. Also for $c\gg 0$, $V^{\al+c+1}M=\sum_it_iV^\al M$ by
the axiom (iii) of \ref{dv}. Therefore the rest of $V^\al M'$ is
recovered from $\sum_{|\nu|\le c}$ through the action of the
$t_i\pa_{u_j}$, hence $V^\al M'$ is finitely generated.

The  axioms (ii) and (iii) of \ref{dv} follow  from the definition
of $V^\al M'$, the simplifications noted above in the
$\DD_{Y'}$-action on $M'$, and the same axioms applied to $V^\al
M$.

The last property to show is the nilpotency of $s-\al$  on
$\Gr_V^\al M'$, where $s=\sum_i\pa_{t_i}t_i+\sum_j\pa_{u_j}u_j$.
Let $m\otimes\pa_u^\nu\in V^\al M'$ with $m\in V^{\al+|\nu|}M$.
Then
$$(s-\al)(m\otimes\pa_u^\nu)=\left(\sum_i\pa_{t_i}t_i-|\nu|-\al\right
)m\otimes\pa_u^\nu.$$
 Hence $(s-\al)^k(m\otimes\pa_u^\nu)\in
V^{\al+1}M'$ if $k$ is the nilpotency order of
$(\sum_i\pa_{t_i}t_i-(\al+|\nu|))$ on $\Gr_V^{\al+|\nu|}M$.
\end{proof}

Examples will be provided in Section \ref{semi}.

\section{Bernstein-Sato polynomials}

The $V$ filtration can be applied to show existence of quite
general Bernstein-Sato polynomials, \cite{BMS}. See \cite{Ka2} for
an account of the classical version of these polynomials.
Following G. Lyubeznik \cite{Ly}, we prove a finiteness result on
the set of all polynomials that are Bernstein-Sato polynomials in
a sense we make precise later.

We keep the notation from the previous section. Suppose first that
$Z\subset X$ is a smooth closed subvariety. To keep this article
as concise as possible we take Theorem \ref{existence} for
granted. Then the quickest way to proceed is by means of the
following technical tool.

\begin{definition}\label{bsu} Let $M$ be  a coherent left $\DD_X$-module. For $u\in M$,
the Bernstein-Sato polynomial $b_u(s)$ of $u$ is the monic minimal
polynomial of the action of $s=-\sum_j\pa_{t_j}t_j$ on
$V^0\DD_Xu/V^1\DD_Xu.$
\end{definition}

We suppressed from the notation the fact that $b_u(s)$ also depends on $Z$. Then we can make explicit the $V$-filtration as follows.

\begin{proposition}\label{bsv} (C. Sabbah \cite{Sab}) If the $V$-filtration along $Z$ exists on $M$,
then $b_u(s)$ exists, it is non-zero for all $u\in M$, and has rational coefficients. Moreover
$$V^\al M=\{\ u\in M\ |\ \al\le c\text{ if } b_u(-c)=0\ \}.$$
\end{proposition}
\begin{proof} Suppose first that $u\in V^\al M$. Recall that
$\sum_j\pa_{t_j}t_j-\beta$ is nilpotent on $V^\beta/V^{>\beta}$
and $V$ is indexed discretely. Then, for a given $\beta$ there is
a polynomial $b(s)$ depending on $\beta$, having all roots $\le -\al$ (and rational), and such that
$b(-\sum_j\pa_{t_j}t_j)\cdot u\in V^\beta$. Hence it is enough to
show that there is $\beta$ such that $V^\beta\cap V^0\DD_Xu\subset
V^1\DD_Xu$.

Let $A=\bigoplus_{i\ge 0}V^i\DD_X\tau^{-i}$ and define
$F_k(A)=\bigoplus_{i\ge 0}(V^i\DD_X\cap \DD_{X,k})\tau^{-i}$. Then
by \ref{noeth}, $A$ is a noetherian ring. Now $\bigoplus_{i\ge
0}V^iM$ is coherent over $A$ because by axiom (iii) of \ref{dv},
there exists $i_0$ such that $V^{i}M$ is recovered from $V^{i_0}M$
if $i\ge i_0$. Denote by $N$ the $V^0\DD_X$-submodule $V^0\DD_Xu$,
and let $U^i=V^i\cap N$ for $i\ge 0$. Then $\bigoplus_{i\ge
0}U^iN$ is also coherent over $A$ since $A$ is noetherian. It
follows that $\bigoplus_{i\ge 0}\Gr_U^iN$ is coherent over
$\bigoplus_{i\ge 0}\Gr_V^i\DD_X$, in particular locally finitely
generated. If $i$ is big compared with the degrees of local
generators, we see that $U^iN\subset V^1\DD_Xu$.

Conversely, fix an element $u\in M$ and suppose that $\al\le c$ whenever $b_u(-c)=0$. Let
$\al_u=\max\{\beta\ |\ u\in V^\beta\}$. We need to show that
$\al\le\al_u$. It is enough to show that $b_u(-\al_u)=0$. For
$\beta\ne\al_u$, $(\sum_j\pa_{t_j}t_j-\beta)$ is invertible on
$V^{\al_u}/V^{>\al_u}$. But $b_u(-\sum_j\pa_{t_j}t_j)u\in
V^{>\al_u}$. Hence we must have $b_u(-\al_u)=0$.
\end{proof}

\begin{lemma}\label{noeth} (\cite{Ka2}-A.29.)
Let $A$ be a be a filtered ring (sheaf on $X$). Assume that
$F_0(A)$ and $\Gr^F(A)$ are noetherian rings, and that
$\Gr_k^F(A)$ are (locally) finitely generated $F_0(A)$-modules for
all $k$. Then $A$ is noetherian.
\end{lemma}

Now let $ Z \subset X $ be an arbitrary closed subset. Let $
f_{1},\dots, f_{r} $ be generators of the ideal of $ Z $, where $
f_{j} \ne 0 $ for any $ j $. Then $ \DD_{X} $ acts naturally on $
\OO_{X}[\prod_{i}f_{i}^{-1},s_{1},\dots, s_{r}]
\prod_{i}f_{i}^{{s}_{i}} $, where the $ s_{i} $ are independent
variables. Define a $ \DD_{X} $-linear action of $ t_{i} $ by $
t_{i}(s_{j}) = s_{j} + 1 $ if $ i = j $, and $ t_{i}(s_{j}) =
s_{j} $ otherwise. Let $ s_{ij} = s_{i}t_{i}^{-1}t_{j} $, and $ s
= \sum_{i}s_{i} $. We will see in Lemma \ref{bu=bz} that under a
well-defined isomorphism the $t_i$'s here correspond to the
$t_i$'s introduced in the second part of section 1.

\begin{definition}\label{defbs}(\cite{BMS}) The Bernstein-Sato polynomial  $
b_{f}(s) $ of $ f := (f_{1}, \dots, f_{r}) $ is defined to be the
monic polynomial of the lowest degree in $ s $ satisfying the
relation
\begin{equation}\label{bs}
b_{f}(s)\prod_{i}f_{i}^{{s}_{i}} = \sum_{j} (P_{j}f_{j}
\prod_{i}f_{i}^{{s}_{i}}),
\end{equation}
where the $ P_{j} $ belong to the ring generated by $ \DD_{X} $
and the $ s_{ij} $. For $ h \in \OO_{X} $,  define similarly $
b_{f,h}(s) $ with $ \prod_{i}f_{i}^{{s}_{i}} $ replaced by $
\prod_{i}f_{i}^{{s}_{i}}h $.
\end{definition}

\noindent {\bf Examples.}

(i) $f=x^2+y^3$. Then $b_f(s)=(s+1)(s+5/6)(s+7/6)$ and
$$P=(\pa^3_y/27+y\pa^2_x\pa_y/6+x\pa^3_x/8).$$

(ii) $f=(x_2x_3,x_1x_3,x_1x_2)$. Then $b_f(s)=(s+3/2)(s+2)^2$ and
the $s_{ij}$ cannot be avoided by the operators $P_j$ in the above
definition (see \cite{BMS}).

\medskip

The polynomial $ b_{Z}(s) := b_{f}(s-r) $ with $ r = \codim_{X}Z $
is shown in \cite{BMS} to depend only on $Z$ and not on $f$. The
existence of non-zero $ b_{f,h}(s) $ follows from the following.

\begin{lemma}\label{bu=bz} With the notation as in \ref{graph},
 if $M=i_*\OO_X$, $u=h\otimes 1$ with $h\in\OO_X$, and the $V$-filtration is taken along $X\times\{0\}$, then
$b_u(s)=b_{f,h}(s)$.
\end{lemma}
\begin{proof}
It suffices to show that $b_u(s)$ is the minimal polynomial of the
action of $s=\sum_js_j$ on
$$\DD_X[s_{ij}]\prod_j
f_j^{s_j}h/\sum_k\DD_X[s_{ij}]f_k\prod_j f_j^{s_j}h,$$ a quotient
of submodules of $\OO_{X}[\prod_{i}f_{i}^{-1},s_{1},\dots, s_{r}]
\prod_{i}f_{i}^{{s}_{i}}h$. We can check that
$\DD_X[s_{ij}]\prod_j f_j^{s_j}h$ and
$\sum_k\DD_X[s_{ij}]f_k\prod_j f_j^{s_j}h$ are isomorphic to
$V^0\DD_Yu$ and $V^1\DD_Yu$. The action of $t_j$ is defined by
$s_j\mapsto s_j+1$, $\prod_j f_j^{s_j}h$ corresponds to $u$, $s_j$
corresponds to $-\pa_{t_j}t_j$, and $s_{ij}=s_it_i^{-1}t_j$.
\end{proof}

Hence it follows from \ref{bsv} that $b_{f,h}(s)$ are
polynomials with rational coefficients. One is allowed to change
the field of definition in the coefficients of the $f_j$'s.

\begin{proposition} There exist non-zero Bernstein-Sato polynomials
$b_{f,h}(s)$ even if in \ref{defbs} one replaces  $X$, $\OO_X$,
and $\DD_X$ with $\cA^n_k$, $k[x_1,\ldots, x_n]$, and $A_n(k)$
(the Weyl algebra) respectively, for $k$ a field of characteristic
zero.
\end{proposition}

\begin{proof} First, suppose that the coefficients of the $f_j$'s
lie in a subfield $K$ of $\CC$. Then also the scalar coefficients
of the $P_j$'s can be assumed to lie in $K$. Indeed, (\ref{bs})
implies, after equating coefficients of monomials in $s_i$'s and
$x_i$'s, that certain $K$-linear relations (*) hold among the
scalar coefficients of the $P_j$'s. Let $L$ be the field generated
by the coefficients of the $P_j$'s. Fix a basis $S$ of $L/K$
containing $1$ and such that every scalar coefficient $c$ which
appears in a $P_j$ can be written as a unique $K$-linear
combination of a finite number of elements of $S$. Let $c_1\in K$
be the coefficient of $1$ in $c$ under this basis, and let
$P_{j,1}$ be the induced operator. Then the $K$-linear relations
(*) hold with $c_1$ replacing $c$, and so (\ref{bs}) holds with
$P_{j,1}$ replacing $P_j$.

Now, going back to our proposition, the conclusion follows from
the Lefschetz principle. Indeed, let $K$ subfield of $k$ generated
over $\QQ$ by the coefficients of the $f_j$'s. Since $\CC$ has
infinite transcendental dimension over $\QQ$, $K$ can be embedded
into $\CC$. Then the coefficients of the $P_j$ are in $K\subset
k$.

\end{proof}

We extend a result of G. Lyubeznik \cite{Ly} to the case of these
more general Bernstein-Sato polynomials. The proof follows closely
his proof.

\begin{proposition}
Fix $n$ and $d$ positive integers. The set of all polynomials
which are of the form $b_f(s)$ for some $f=f_1,\ldots, f_r\in
k[x_1,\ldots, x_n]$ with $\deg f_i\le d$ is finite even if $k$ is
varying over all the fields of characteristic zero.
\end{proposition}

\begin{proof}. Let $N$ be the
number of monomials in $x_1,\ldots, x_n$ of degree $\le d$. Then
the $f$'s are the closed $k$-rational points of
$[\mathbf{A}^N_k]^{\times r}$. Let $P=\sum_{|\alpha |\le
d}c_\alpha x^\alpha$ be the polynomial of $n$ variables of degree
$d$ with underterminate coefficients. Then $B_k=k[\text{ the
}c_\alpha{\text 's}]^{\otimes r}$ is the coordinate ring of
$[\mathbf{A}^N_k]^{\times r}$. Define $$F_i=1\otimes\ldots\otimes
P\otimes\ldots\otimes 1 \in B_{\mathbf{Q}}[x_1,\ldots, x_n]$$ by
placing $P$ in the $i$-th position. Here $\times$ and $\otimes$
mean over $\mathbf{Q}$.

Let $Y$ be a reduced and irreducible closed subset of
$[\mathbf{A}^N_\mathbf{Q}]^{\times r}$. Denote by $G=(G_i)_i$ the
image of $F=(F_i)_i$ under the natural $\mathbf{Q}$-algebra
homomorphism $B_{\mathbf{Q}}[x_1,\ldots, x_n]\rightarrow
\mathbf{Q}(Y)[x_1,\ldots,x_n]$, where $\mathbf{Q}(Y)$ is the
function field of $Y$. Let $\mathbf{Q}[Y]$ be the coordinate ring
of $Y$. Then,  we have a functional equation

\begin{equation}\label{fe}
b_G(s)\prod_iG_i^{s_i}=\sum_j P_jG_j\prod_iG_i^{s_i}
\end{equation}
with $s=\sum s_i$ and $P_j\in A_n(\mathbf{Q}(Y))[\text{ the
}s_{ij}\text{'s}]$. Denote by $c\in \mathbf{Q}[Y]$ the common
denominator. Denote by $U(c)\subset Y$ the subscheme whose
coordinate ring is $\mathbf{Q}[Y]_c$. By specializing (\ref{fe}),
the $f$'s given by the closed $k$-rational points of $U(c)\times
k$ have $b$-functions $b_f(s)$ dividing $b_G(s)$. Hence they are
only finitely many such $b_f(s)$ even if $k$ varies.

We proceed now by induction on the dimension of $Y$ proving that
the $k$-rational points of $Y\times k$ give only finitely many
$b$-functions even if $k$ varies. For dimension zero, $c=1$ and so
$U(c)=Y$. In higher dimensions, $Y-U(c)$ is the union of reduced
and irreducible closed subsets of smaller dimension.
\end{proof}

\section{Multiplier ideals}\label{semi}

The multiplier ideals introduced by A. Nadel \cite{Na} encode the
complexity of singularities via their resolutions. It turns out
that they are essentially the same as the $V$ filtration on
$\OO_X$.

Let $X$ be a smooth complex variety and $Z\ne X$ a closed
subscheme. Let $\mu:X'\ra X$ be a log resolution of $(X,Z)$. That
is $\mu$ is proper birational, $X'$ is smooth, and
$\Ex(\mu)\cup\mu^{-1}Z$ is a divisor with simple normal crossings.
Here $\Ex(\mu)$ denotes the exceptional locus of $\mu$. Let
$\II\subset\OO_X$ denote the ideal of $Z$. Let $H$ be the
effective divisor on $X'$ such that
$\mu^{-1}(\II)\cdot\OO_{X'}=\OO_{X'}(-H)$.

\begin{definition} For  $\al>0$, the multiplier ideal of $(X,\al\cdot Z)$ is
defined as
$$\MI{\al\cdot
Z}=\mu_*(\omega_{X'/X}\otimes\OO_{X'}(-\rndown{\al\cdot
H})).$$
\end{definition}
Here
$\omega_{X'/X}=\det\Omega_{X'}^1\otimes\mu^*(\det\Omega_X^1)^\vee$
is the sheaf of relative top-dimensional forms, and $\rndown{.}$
rounds down the coefficients  of  the irreducible divisors. One
can extend this definition to a formal combination of closed
subschemes $\sum_i\al_i\cdot Z_i$ by replacing $\al\cdot H$ with
$\sum_i\al_i\cdot H_i$.

The original analytic definition of multiplier ideals is, locally,
$$\MI{\al\cdot
Z}=\{\ h\in\OO_X\ |\ |h|^2/(\sum_{1\le i\le r}|f_i|^{2})^\al \in\
L^1_{loc}\ \},$$ where $f_1,\ldots, f_r$ generate $\II$. The first
definition shows the second is independent of the choice of
generators, and the second definition shows the first is
independent of the choice of resolution. See \cite{La} for more on
multiplier ideals.

The multiplier ideals measure how singular $Z$ is. The intuition
here is that smaller multiplier ideals means worse singularities.
For example, varying the coefficient in front of $Z$, one obtains
a decreasing family $\{\MI{\al\cdot Z}\}_{\al\in\QQ}$. Because of
the rounding-down of coefficients in the construction of
$\MI{\al\cdot Z}$ there exist positive rational numbers $ 0 <
\alpha_{1} < \alpha_{2} < \cdots $ such that $ \MI{\alpha_{j}\cdot
Z} = \MI{\alpha \cdot Z} \ne \MI{\alpha_{j+1}\cdot Z} $ for $
\alpha_{j} \le \alpha < \alpha_{j+1} $ where $ \alpha_{0} = 0 $.
These numbers $ \alpha_{j} \,(j > 0) $ are called the jumping
numbers of the multiplier ideals associated to $ (X,Z) $. The
log-canonical threshold of $(X,Z)$ is the smallest non-zero
jumping number, $\lc(X,Z)=\al_1.$ Equivalently, $\lc(X,Z)$ is the
number $\al$ such that $\MI{\al\cdot Z}\ne\OO_X$, but
$\MI{(\al-\eps)Z}=\OO_X$ for $0<\eps\ll 1$.

\medskip

\noindent {\bf Examples.}

(i) $Z=\{\ x^2+y^3=0\ \}\subset\cA^2$. Then $\MI{\al\cdot Z}$ is
equal to $\OO_X$ if $0<\al<5/6$, and it is the maximal ideal at
$(0,0)$ if $5/6\le\al<1$.

(ii) $Z=\{\ x_1x_2=x_2x_3=x_1x_3=0\ \}\subset\cA^3$. Then
$\MI{\al\cdot Z}$ is equal to $\OO_X$ if $0<\al<3/2$, and it is
the ideal $(x_1,x_2,x_3)$ if $3/2\le\al <2$. This follows from
\cite{La}-III.9.3.4 which gives the formula for multiplier ideals
of monomial ideals.

(iii) The multiplier ideals of hyperplane arrangements, and more
generally, stratified locally conical divisors are determined in
\cite{Mu}, and respectively \cite{S04}.

\begin{theorem}\label{miv}(\cite{BMS}, \cite{BS}.) For $\al>0$,
$V^\al\OO_X=\MI{(\al-\eps)\cdot Z}$, where $0<\eps\ll 1$ and the
filtration $V$ of $\OO_X$ is taken along $Z$ as in \ref{graph}.
\end{theorem}

The relation with Bernstein-Sato polynomials is then given by
Proposition \ref{bsv} and Lemma \ref{bu=bz}:

\begin{corollary}\label{mibs} For $\al >0$, $$\MI{\al\cdot Z}=\{
h\in\OO_X\ |\ \al<c\text{ if } b_{f,h}(-c)=0\ \},$$ where
$f=f_1,\ldots, f_r$ is any set of generators of the ideal
$\II\subset\OO_X$ of $Z$.
\end{corollary}

In particular, $$\lc(X,Z)=-(\text{biggest root of }b_f(s)),$$
since $b_f(s)=b_{f,1}(s)$ by definition \cite{Ko}, \cite{Li}.

\section{Monodromy of vanishing cycles}

The initial scope of the $V$ filtration of M. Kashiwara and B.
Malgrange was to to construct vanishing cycles in the category of
(regular holonomic) $\DD$-modules.

Let $X$ be a smooth complex variety. Denote by $M_{rh}(\DD_X)$ the
abelian category of regular holonomic $\DD_X$-modules, and by
$D^b_{rh}(\DD_X)$ the derived category of bounded complexes of
$\DD_X$-modules with regular holonomic cohomology. By A.
Beilinson, this is equivalent with the bounded derived category of
$M_{rh}(\DD_X)$. Let $D^b_c(X)$ be the derived category of bounded
complexes of sheaves (in the analytic topology of $X$) of
$\CC$-vector spaces with constructible cohomology. The
Riemann-Hilbert correspondence generalizing the analogy between
the $\DD_X$-module $\OO_X$ and the constant sheaf $\CC_X$ states
(see \cite{Bo}):

\begin{theorem}(M. Kashiwara, Z. Mebkhout) Let $X$ be a smooth complex variety.
  There is a well-defined functor
$$\te{DR}:D^b_{rh}(\DD_X)\lra D^b_c(X)$$ which is an equivalence of
categories commuting with the usual six functors. $\te{DR}$ also
defines an equivalence $M_{rh}(\DD_X)\ra\text{Perv}\,(X)$, where
$\text{Perv}\,(X)\subset D^b_c(X)$ is the subcategory of perverse
sheaves.
\end{theorem}

Let $f\in\OO_X$ be a regular function. The vanishing cycles
functor $\phi_f$ on $D^b_c(X)$ and the monodromy action $T$ on it
should then have a meaning only in terms of $\DD$-modules since
the shift $\phi_f[-1]$ restricts as a functor to $\text{Perv}\,
(X)$.  Let $M$ be a regular holonomic $\DD_X$-module. Let $\ti{M}$
be its direct image under the graph of $f$, as in section 2. If
$M$ is also quasi-unipotent then there exists a $V$ filtration
indexed by $\QQ$ on $\ti{M}$ along $X\times\{0\}$. If $M$ is not
quasi-unipotent, a close version of the following still holds:

\begin{theorem}\label{van}(M. Kashiwara, B.
Malgrange ) Let $\al\in [0,1)$ be a rational number.
$\Gr_V^\al\ti{M}$ corresponds to the $\exp (-2\pi
i\al)$-eigenspace of $\phi_f[-1](\te{DR}(M))$ with respect to the
action of the semisimple part $T_s$ of the monodromy.
\end{theorem}

Combining this result with Theorem \ref{miv} and with additional
structures such as mixed Hodge modules, one obtains a relation
between multiplier ideals and the Hodge spectrum of hypersurface
singularities. Let $f:X\ra \cA^1$ be a regular function. Recall
that if $i_x:x\hookrightarrow f^{-1}(0)$ is a point, the Milnor
fiber of $f$ at $x$ is is the Milnor fiber of the corresponding
holomorphic germ $f:(\CC^m,0)\ra (\CC,0)$
$$M_{f,x}=\{z\in\CC^{m}\ |\ |z|<\eps\ {\rm and\ }f(z)=t\}$$ for a
fixed $t$ with $0<|t|<\eps\ll 1$. Then
\begin{align}\label{mf}
H^i(i_x^*\phi_f\CC_X)&=\ti{H}^{i}(M_{f,x},\CC),
\end{align}
where $\ti{H}$ stands for reduced cohomology. These vector spaces
are endowed with the monodromy action $T$ and with mixed Hodge
structures on which $T_s$ acts as automorphism. Indeed, the mixed
Hodge module theory of M. Saito on the left-hand side of
(\ref{mf}) recovers the mixed Hodge structure of V. Navarro-Aznar
from the right-hand side. As numerical invariants encoding the
behaviour of the Hodge filtration $F$ under $T_s$ one has the
generalized equivariant Euler characteristics
$$n(i,\al)=\sum_j(-1)^j\dim\Gr_F^{i}\ti{H}^j(M_{f,x},\CC)_\al,$$
where $\al\in \QQ\cap [0,1)$, $i\in\{0,\ldots, m-1\}$, and the
subscript $\al$ stands for the eigenspace of $T_s$ with eigenvalue
$exp(2\pi i\al)$. These invariants form the Hodge spectrum of $f$
introduced by J. Steenbrink \cite{St}. For $\al\in (0,1]$, let
$$n_{\al, x}(f)=(-1)^{n-1}n(m-1,1-\al),$$
so that $n_{\al, x}(f)$ describe the spectrum for the smallest
piece of the Hodge filtration.

\medskip
\noindent {\bf Example.} If $f=x^2+y^3$ and $x=(0,0)\in \cA^2$,
then $n_{\al,x}(f)$ is zero for $\al\notin \{5/6, 7/6\}$, and is 1
otherwise.

\medskip

On the other hand, for every jumping number $\al\in (0,1]$ of
$(X,Z)$ where $Z$ is the zero set of a regular function $f$,
define the inner jumping multiplicity at $x$
$$n_{\al,x}(Z)=\dim \MI{{(1-\eps)\al\cdot Z
}}/\MI{{(1-\eps)\al\cdot Z}+\delta\cdot x},$$ where
$0<\eps\ll\delta\ll 1$. It is proved in \cite{B} that
$n_{\al,x}(Z)$ is finite and does not depend on $\eps$ and
$\delta$. Let $\widetilde{\OO_X}$ be  the $\DD$-module direct
image of $\OO_X$ under the graph of $f$. In connection with
\ref{miv} it is crucial to observe that the smallest piece of the
Hodge filtration on $V^\al\widetilde{\OO_X}$is exactly
$V^\al\OO_X$. Then the above arguments lead to:

\begin{theorem}( \cite{B}, \cite{BS} ) For $\al\in (0,1]$,
$$n_{\al,x}(f)=n_{\al,x}(Z).$$
\end{theorem}

\bibliographystyle{plain}
\bibliography{survey}

\end{document}